\documentclass{article}
\usepackage[utf8]{inputenc}
\usepackage{eurosym}
\usepackage{amstext} 
\DeclareRobustCommand{\officialeuro}{%
  \ifmmode\expandafter\text\fi
  {\fontencoding{U}\fontfamily{eurosym}\selectfont e}}
\usepackage{nomencl}
\usepackage{etoolbox}
\usepackage{amsmath}
\usepackage{amssymb}
\usepackage{array}
\usepackage{mdwmath}
\usepackage{mdwtab}
\usepackage{array}
\usepackage{graphicx}
\usepackage{appendix}
\usepackage{cleveref}
\usepackage{textcomp}
\usepackage{caption}
\usepackage{subcaption}
\usepackage{comment}
\usepackage[left=3cm, right=3cm, top=2.5cm]{geometry}

\newcommand{\olsi}[1]{\,\overline{\!{#1}}} 

\renewcommand\nomgroup[1]{%
	\ifstrequal{#1}{I}{\item[\textbf{Indices and Sets}]}{%
	\ifstrequal{#1}{P}{\vspace{10pt} \item[\textbf{Parameters}]}{%
	\ifstrequal{#1}{V}{\vspace{10pt} \item[\textbf{Variables}]}{}}}%
	\ifstrequal{#1}{W}{\vspace{10pt} \item[\textbf{Dual variables}]}{}
}

\title{Grid Tariffs for Peak Demand Reduction: Is there a Price Signal Conflict with Electricity Spot Prices?}

\author{Sigurd Bjarghov$^{a*}$, Matthias Hofmann$^{ab}$ \\
$^a$Department of Electric Power Engineering, NTNU - Trondheim, Norway\\
$^b$Statnett - Oslo, Norway\\
$^{*}$Corresponding author - Email:\ sigurd.bjarghov@ntnu.no}

\date{}

\begin{document}




\maketitle

\begin{abstract}

The electricity grid is expected to require vast investments due to the decarbonization-by-electrification trend, calling for a change in grid tariff design which provides proper incentives for reducing peak loads. However, price signals from grid tariffs could be “distorted” from electricity spot prices which also represents a significant of the total consumer electricity bill. This paper attempts to identify whether there is a price signal conflict between grid tariffs and spot prices. Four different grid tariff designs are compared, using a generic demand response model as part of a cost-minimizing linear program to simulate the reduction in peak load. The method is applied to metered electricity demand from 3608 consumers in Oslo, Norway. Results show that new grid tariff designs reduce peak loads by 1-4\%, and that reduction in peak load is smaller when consumers are subject to electricity spot prices.

\end{abstract}


\makenomenclature
\setlength{\nomlabelwidth}{1.5cm}
\nomenclature[I]{$T$}{Set of hourly time steps, index $t$}
\nomenclature[I]{$D$}{Set of days, index $d$}
\nomenclature[I]{$M$}{Set of months, index $m$}

\nomenclature[P]{$L_{t}$}{Original load series [kWh/h]}
\nomenclature[P]{$C^{red}$}{Discomfort cost of reducing load [$\frac{ \euro }{kWh}$]}

\nomenclature[P]{$C^{ET}$}{Energy term [$\frac{ \euro }{kWh}$]}
\nomenclature[P]{$C^{h}$}{Energy term above subscribed capacity [$\frac{ \euro }{kWh}$]}
\nomenclature[P]{$C^{sub}$}{Annual subscribed capacity cost[$\frac{\euro }{kW \times year}$]}
\nomenclature[P]{$C^{spot}$}{Electricity spot price [$\frac{ \euro }{kWh}$]}
\nomenclature[P]{$C_t^{TOU}$}{Static time-of-use network tariff [$\frac{ \euro }{kWh}$]}
\nomenclature[P]{$C_t^{dynTOU}$}{Dynamic time-of-use network tariff [$\frac{ \euro }{kWh}$]}
\nomenclature[P]{$C^{peak}$}{Demand charge peak cost[$\frac{ \euro }{kWh}$]}
\nomenclature[P]{$C^{e}$}{Total electricity cost [\euro]}
\nomenclature[P]{$C^{g}$}{Total grid tariff cost [\euro]}
\nomenclature[P]{$C^{f}$}{Total flexibility use cost [\euro]}

\nomenclature[P]{${Q}^{flex}$}{Share of reducible load in an hour [\%]}
\nomenclature[P]{${E}^{flex}$}{Share of reducible energy in a day [\%]}

\nomenclature[V]{$x_{t}^{i}$}{New load series with demand response [kWh/h]}
\nomenclature[V]{$x_{t}^{l}$}{Electricity consumption below sub. cap. [kWh/h]}
\nomenclature[V]{$x_{t}^{h}$}{Electricity consumption above sub. cap. [kWh/h]}
\nomenclature[V]{$x_{m}^{max}$}{Monthly peak demand [kWh/h]}
\nomenclature[V]{$x^{sub}$}{Subscribed capacity [kW]}
\nomenclature[V]{$q_{t}^{red}$}{Load reduction [kWh/h]}



\section{Introduction} \label{sec:introduction}

Aiming to reduce 55\% of carbon emissions by 2030 \cite{klimamal2030}, Norway plans extensive electrification, especially in the transport and industrial sectors \cite{NVE2020,Statnett2020}. As a consequence, the transmission and distribution grid operators are facing vast amounts of connection requests from commercial consumers. As it is demanding to build grids at sufficient speed to incorporate these new grid connection requests, it is of increasing importance to ensure efficient use of existing grids. Meanwhile, the distribution grid is expecting increased peak demand due to electrification of transport and other power-intensive loads which in addition to urbanization trends results in an expected increase in congestions, also in the low voltage grid.

A solution to this development could be achieved by moving from flat, volumetric grid tariff designs to time-of-use or capacity-based grid tariffs. With redesigned grid tariffs, electricity peak loads can be reduced by implicit flexibility, i.e. consumers reacting to price signals by reducing or shifting demand, often referred to as demand response.

The future power system requires more precise pricing mechanisms as integration of demand side flexibility and distributed generation introduces new challenges, especially in the distribution grid. Volumetric grid tariffs do not sufficiently represent cost-reflectivity and are already responsible for inefficient investments and operational decisions \cite{PerezArriaga2016}. Implicit flexibility through price signals also avoids struggles with market manipulation and baseline related issues such as local flexibility markets may \cite{Ziras2021WhyMarkets}.

A variety of grid tariff designs to incentivize peak demand reduction from flexibility have been suggested in recent literature. Optimal time-of-use tariff design is one option, and should preferably be demand-based to achieve peak load reductions \cite{Li2016AModel}. The design of time-of-use tariffs can be difficult, as the welfare increase is dependent on which technologies exist on consumers side in terms of demand response cost \cite{Yang2013ElectricityConsideration}. 

Of the capacity-based grid tariffs, demand charges have received significant attention due to the welfare gains under higher shares of renewable generation in the distribution grid \cite{Brown2018OnResources}. Still, coincidence-related issues, i.e. the lack of guarantee that residential peak and system peak coincide reduces the welfare gain significantly \cite{Borenstein2016}, suggesting that dynamic tariffs which adapt to the grid status are more likely to increase welfare \cite{Hledik2014}. Subscribed capacity tariffs also provide incentives reduce peak loads in neighborhoods with residential consumers and similar local energy systems, both under cooperative \cite{Backe2020} and competitive market conditions \cite{Bjarghov2020}. The grid tariff design not only impacts operation, but also investments in decarbonized neighborhoods \cite{Pinel2019}.

However, grid tariffs are not the only price signal consumers are exposed to. Grid tariffs have historically made up roughly one third of the total electricity bill in Norway, with taxes and electricity spot prices also taking one third each. During the end of 2021 and winter of 2022, Europe has experienced historically high electricity spot prices at a size which easily could "outperform" the most suggested grid tariff structures in the sense that consumers would react to the price signal from the electricity spot prices, rather than from the grid tariff. Analysis on demand response from a combined spot price and grid tariff signal has been proposed in the literature, but are often complex and difficult to implement \cite{Schreiber2015FlexibleGrid}. Consumers on fixed price contracts will not respond to spot prices, but due to historically low prices, more than 95\% of Norwegian residential consumers are on spot price or variable price contracts \cite{SSBcontracttypes}.

This raises the questions: Is there a price signal conflict between electricity spot prices and different grid tariff designs, and how large is it? This is of particular interest as there is often (but not always) a correlation between high electricity prices and cold winters with high demand, which also is the dimensioning factor for grid expansion. In other words, if the cold, premise-setting winters for grid expansions might include very high spot prices, which grid tariff designs are the most robust in order to achieve peak load reduction in those few days which might occur as seldom as every decade?

This article attempts to answer these questions by simulating demand response for peak demand reduction, using historical spot prices and real, metered data from 3608 consumers in Oslo, Norway, from November 2020 to October 2021. Summarized, the main contributions of this article are the following:
\begin{itemize}
    \item A quantification of the price signal conflict between electricity spot prices and grid tariffs, with respect to reducing peak loads.
    \item A comparison of peak demand reduction under different grid tariff designs, when exposed to both real-time electricity spot prices and no spot prices.
\end{itemize}

The remainder of the paper is organized as follows: \Cref{sec:gridtariffs} discusses the different grid tariff designs. The method and optimization model is presented in \Cref{sec:method}, followed by the case study description in \Cref{sec:casestudy}. Results and discussions are then presented in \Cref{sec:results}, followed by conclusions and further work suggestions in \Cref{sec:conclusion}.


\section{Grid tariffs} \label{sec:gridtariffs}

Grid tariffs represent the cost of transferring electricity from the place of generation to consumption. Ideally, they should reflect both the long- and short-term marginal cost of transferring electricity, but this is a burdensome task as the real cost depends on complex mechanisms. Grid tariffs are designed not only to be cost-reflective, but also after a number of a criteria such as cost-efficiency, cost-recovery, complexity, implementation burden, acceptance, transparency and fairness \cite{ACER2021,EDSO2021FutureGuidance}. Often, there is a reverse relationship between these criteria, making tariff design a task of finding the least-worse alternative with respect to all the criteria \cite{PerezArriaga2016}. 

Traditionally, grid tariffs have been designed to transfer costs from the distribution system operators to the consumers in a simple manner. Volumetric tariffs have done this job relatively well with respect to simplicity and cost-recovery, but is especially inefficient in terms of incentivizing flexibility response from consumers for peak reduction, which is the main focus of this article.

Grid tariffs are often split into three types of costs: a) fixed, b) volumetric and c) capacity-based costs, aiming to represent different types of costs related to administration, as well as short and long-term marginal costs of electricity consumption \cite{Hennig2022WhatAssessment}. In this article, we look into two time-of-use tariffs, as well as two types of capacity-based tariffs. The list of tariff models is based on literature and proposed tariffs in Norway, limited to distribution grid tariffs. The included tariffs are presented below, whereas the cost levels can be found in \Cref{tab:GTcosts}. 

\subsection*{Subscribed capacity}
Capacity-subscription tariffs are based on consumers subscribing to a capacity ex-ante, which has a cost per kilowatt, for example annually. Consumption below the subscribed level is subject to a small energy term, often reflective of the marginal losses in the grid, whereas consumption above the level is subject to an excess energy term which is significantly higher. Consumers then have an incentive to stay below the subscribed capacity.

\subsection*{Demand charges}
Demand charges are based on the consumers peak demand in a given time period, typically monthly. Consumers have then an incentive to have a low peak demand, which is measured as the highest electricity use in one hour.

\subsection*{Static time-of-use}
Static time-of-use has a volumetric cost part with predetermined energy cost that can shift from hour to hour, aiming to incentivize use when the demand is low and similarly penalizing consumption when the demand is high and possible congestions in the grid occur. Typically, this involves having a higher price per kilowatt-hour in the morning and in the evening, with the option of seasonal variation.

\subsection*{Dynamic time-of-use}
Unlike the static version, the dynamic version is only active when there is scarcity in the grid. This can be defined as a certain number of days per year with the highest grid utilization, adding a very high energy term during the peak load hours of those days.


\section{Method}\label{sec:method}

\subsection{Approach} \label{subsec:methodframework}

The method used in this paper is to simulate the total demand response from a large set of consumers. The demand response model is described in \Cref{subsec:optimizationmodel}. The simulation is performed by formulating an optimization model which minimizes costs of each individual consumer, formulated as a linear program. This results in a new demand curve after consumers have tried to reduce costs using the modeled demand response, which then is used to discuss the efficiency of the different grid tariff designs, both with and without being subject to spot prices. The optimization model is introduced in \Cref{subsec:optimizationmodel}, whereas the specific data and simulated cases are presented in \Cref{sec:casestudy}.

\subsection{Optimization model} \label{subsec:optimizationmodel}

The consumer problem is formulated as a linear cost minimizing program aiming to minimize the sum of electricity costs $C^e$ (if applicable), grid tariff costs $C^g$ and flexibility usage costs $C^f$ as shown in \eqref{eq:objective}. The full nomenclature can be found in the Appendix.

\begin{equation}
    \min C^{e} + C^{g} + C^{f} \label{eq:objective}
\end{equation}

The total electricity cost $C^e$ is given by \eqref{eq:totalelectricitycost}, whereas the flexibility cost $C^f$ is described in \eqref{eq:curtflexcost}. 

\begin{equation}
    C^{e} = \sum_t x_{t}^{i} C_{t}^{spot} \label{eq:totalelectricitycost}
\end{equation}

The grid tariff costs $C^g$ are described in the following subsections.

\subsubsection*{Subscribed capacity tariff} \label{subsec:}

Capacity subscription tariffs involve consumers taking an active choice where they subscribe to a capacity level ($X^{sub}$) once a year, which is associated with a capacity cost $C^{sub}$ per kilowatt. Electricity consumption of the end-user is split into demand below ($x_{t}^{l}$) and above ($x_{t}^{h}$) the subscription level as shown in in \eqref{subcap1} and \eqref{subcap2}. Under this grid tariff structure, we assume that all consumers subscribe to their optimal level which is found by setting $X^{sub}$ as a variable in the consumer problem. However, this level will not be the same with and without flexibility assets. We therefore find the optimal subscribed capacity first, and then set this value as a fixed parameter in the problem again when flexibility is added.

\begin{subequations}
\begin{align}
    & x_{t}^{l} + x_{t}^{h} = x_{t}^{i} c \label{subcap1} \\
    & x_{t}^{l} \leq X^{sub} \quad \forall t \label{subcap2}
\end{align}
\end{subequations}

Finally, the grid cost function is given by \eqref{objfunsubcap}, where demand below and above the subscribed capacity are associated with the energy cost term $C^{ET}$ and excess energy cost term $C^h$, respectively. 

\begin{equation}
    C^{g} = x^{sub}\cdot C^{sub} + \sum_{t}  (x_{t}^{l}C^{ET} + x_{t}^{h}C^{h})  \label{objfunsubcap} 
\end{equation}

\subsubsection*{Demand charge tariff}

Demand charges penalizes the monthly peak demand of the consumer $x_{m}^{max}$ by a specific cost per kW peak $C^{peak}$. The cost function is given by \eqref{objfunMP}. The general energy term $C^{ET}$ is added, similar to the other tariffs. Additionally, another constraint to enforce the peak demand cost is needed as shown in \eqref{eq:peakdemandcost}.

\begin{equation}
    C^{g} = \sum_{m} x_{m}^{max} C^{peak} + \sum_{t} x_t C^{ET} \label{objfunMP} 
\end{equation}

\begin{equation}
    x_{t}^{i} \leq x_{m}^{max} \quad \forall t \label{eq:peakdemandcost}
\end{equation}

\subsubsection*{Time-of-use tariff}

The time-of-use tariff adds a specific cost of using electricity at different time steps as shown in \Cref{tab:GTcosts}. Outside the peak price hours, there is a volumetric energy term $C^{ET}$. 

\begin{equation}
    C^{g} = \sum_{t} x_{t}^{i}C_t^{TOU} \label{objfunTOU} 
\end{equation}

\subsubsection*{Dynamic time-of-use tariff}

The dynamic time-of-use tariff has a very high cost term for a selected number of days per year, based on which days have the highest peak loads. Outside those hours, the price has a regular energy term $C^{ET}$.

\begin{equation}
    C^{g} = \sum_{t} x_{t}^{i}C_t^{dynTOU} \label{objfunTOUdyn} 
\end{equation}

\subsection{Flexibility model}

The cost of flexibility $C^{f}$ is given by \Cref{eq:curtflexcost}, which adds a cost $C^{red}$ per kilowatt-hour of reduced electricity demand. This is not a monetary cost, but a discomfort cost, representing the discomfort of responding to price signals.

\begin{align}
    & C^{f} = \sum_{t} q_{t}^{red}\cdot C^{red} \label{eq:curtflexcost}
\end{align}

The energy balance is given by \eqref{energybalancecurt}, where the new load series $x_{t}^{i}$ is the sum of the original load $L_{t}$ minus the demand reduction $q_{t}^{red}$.

\begin{align}
     & x_{t}^{i} = L_{t} - q_{t}^{red} \quad \forall t \label{energybalancecurt}
\end{align}

The flexibility is modeled as the ability to reduce load without shifting to other hours and is modeled in generic terms rather than as assets. The advantage of this modeling approach is the ability to emulate a general demand response from a set of consumers, without the computational efforts of asset modeling. It also draws advantage from not assuming what kind of flexibility assets that exist, or will exist in the future. Instead, the model represent a generic flexibility response specified by 2 parameters:

\begin{itemize}
  \item Max. possible power reduction in an hour, $Q^{flex}$, in \% of demand in that hour
  \item Max. possible electricity demand reduction in a day, $E^{flex}$, in \% of demand that day
\end{itemize}

These two parameters set a limit to how much power can be reduced in an hour, as well as how much energy can be reduced per day, as shown in \eqref{curtmaxpower} and \eqref{curtmaxenergy}, respectively. The values of these parameters are determined based on results of international studies \cite{VaasaETT2017, Faruqui2017} as well as Norwegian studies \cite{Hofmann2021, Siebenbrunner2022}.

\begin{subequations}
\begin{align}
     & q_t^{red} \leq Q^{flex} \cdot L_t \quad \forall t \label{curtmaxpower} \\
     & \sum_{t\in d} q_t^{red} \leq E^{flex} \cdot \sum_{t\in d} L_t \quad \forall d \label{curtmaxenergy}
\end{align}
\end{subequations}
 It is assumed that such a representation of flexibility is relatively accurate when modeling large sets of consumers, although the spread in flexibility response from each individual consumer obviously is not equal as assumed in this study. Additionally, the following assumptions were made:

\begin{itemize}
  \item The discomfort cost parameter $C^{red}$ is high enough to reflect an assumed discomfort cost of being flexible, but small enough to trigger activation for all tariffs.
  \item If several hours are equally optimal for reduction (ToU, subscribed capacity), the relative reduction is equal in all these hours.
\end{itemize}

\section{Case Study \& Data} \label{sec:casestudy}

The case study and data are similar as in \cite{HofmannPreprint}, where hourly electricity load data from 3608 consumers, including 3081 household consumers and 527 commercial consumers, were used for the analysis. The households stand for 47 \% of the electricity demand, and the commercial consumers for 53 \% respectively.
The consumers are located at 112 different substations, under the same transformer in Oslo, Norway, and the metered data are from the period November 2020 to October 2021. In addition, historical spot prices and demand from the price zone NO1 are collected for the same period. These data period is of particular interest because the winter of 2020/2021 saw the highest measured electricity consumption in Norway.



Several cases are studied to understand the difference in achieved grid peak demand reduction, and are listed below. They include two benchmark tests, as well as a performance analysis of the four grid tariffs described in \Cref{sec:gridtariffs}. The grid tariff case includes the test of the grid tariffs without any spot price signal. This is expanded in the grid tariff and spot price case, which includes the combined price signal from each of the four grid tariffs as well as the spot price. Finally, the system optimal response and the spot price case represent the benchmark tests. The system optimal response case is defined as the system's maximum ability to reduce peak load with the available flexibility, and assuming that reduction is coordinated between all consumers. The spot price case is the cost minimizing response under spot prices only. In summary the following cases are analyzed:
\begin{itemize}
    \item System optimal response case
    \item Grid tariff case
    \item Grid tariffs plus spot price case
    \item Spot price case
\end{itemize}

\begin{table}[h]
\caption{Cost levels for all grid tariffs.}
\begin{tabular}{lll}
\hline
\textbf{Symbol} & \textbf{Cost level} & \textbf{Comment}                                                                                                  \\ \hline
$C^{ET}$            & 0.25 $\frac{NOK}{kWh}$        & Applied to all tariffs                                                                                            \\
$C^{TOU}$           & 1.2 $\frac{NOK}{kWh}$         & \begin{tabular}[c]{@{}l@{}}06-22 during winter except weekends\\ and holidays, energy term otherwise\end{tabular} \\
$C^{dynTOU}$        & 4.5 $\frac{NOK}{kWh}$         & 06-22, during the 20 peak load days                                                                               \\
$C^{sub}$           & 1000 $\frac{NOK}{kW-year}$    & \begin{tabular}[c]{@{}l@{}}Cost per kilowatt subscribed \\ capacity per year\end{tabular}                                             \\
$C^{h}$             & 1.65 $\frac{NOK}{kWh}$        & \begin{tabular}[c]{@{}l@{}}Cost for electricity above subscription \\ level, energy term when under\end{tabular}  \\
$C^{peak}$          & 75 $\frac{NOK}{kW-month}$     & Cost for peak load per month                                                                                      \\ \hline
\end{tabular}
\label{tab:GTcosts}
\end{table}

As the cost-recovery should be similar if there is no demand response, all the cost parameters of grid tariffs are determined using backwards calculation with respect to the existing energy tariff. The cost levels are shown in \Cref{tab:GTcosts}.

The parameters $Q^{flex}$ and $E^{flex}$ for available demand flexibility were chosen based on the aforementioned studies on demand response and are set to 25\% and 2.5\%, respectively.

\section{Results \& discussion}\label{sec:results}

\subsection{Simultaneity of peak demand and peak spot prices}
The spot price is the result of the market equilibrium between supply and demand bids on the day ahead market. Therefore, everything else equal, a higher electricity demand implies higher spot prices. In addition, other factors, as for example fuel prices and weather, change the supply bids and in consequence the spot price. These effects can also be observed in the historical demand and spot prices in the case study.

\begin{figure}
     \centering
     \begin{subfigure}[b]{0.35\textwidth}
         \centering
         \includegraphics[width=0.95\textwidth]{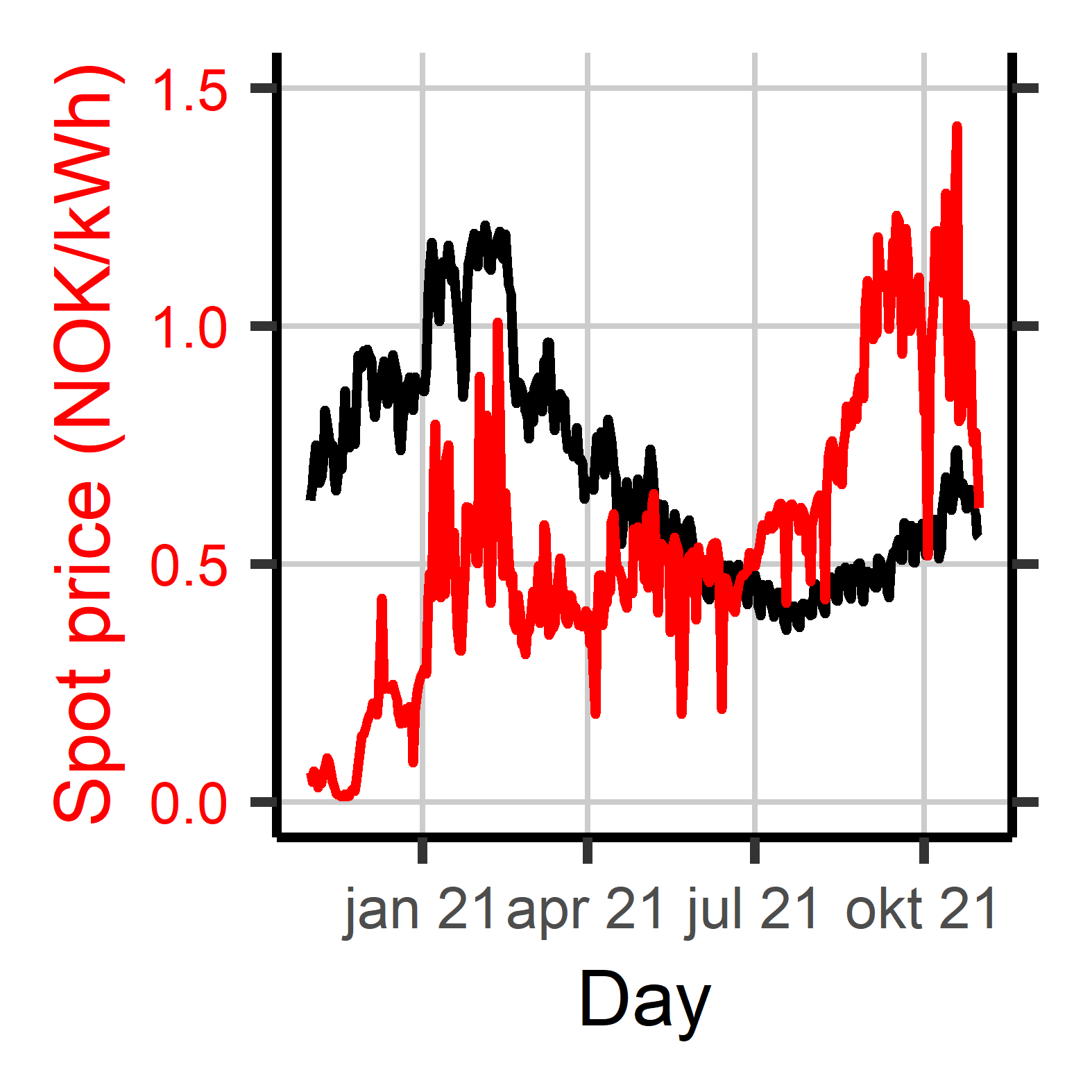}
         \caption{Daily}
         \label{fig:yearlycomparison}
     \end{subfigure}
     \hfill
     \begin{subfigure}[b]{0.35\textwidth}
         \centering
         \includegraphics[width=0.95\textwidth]{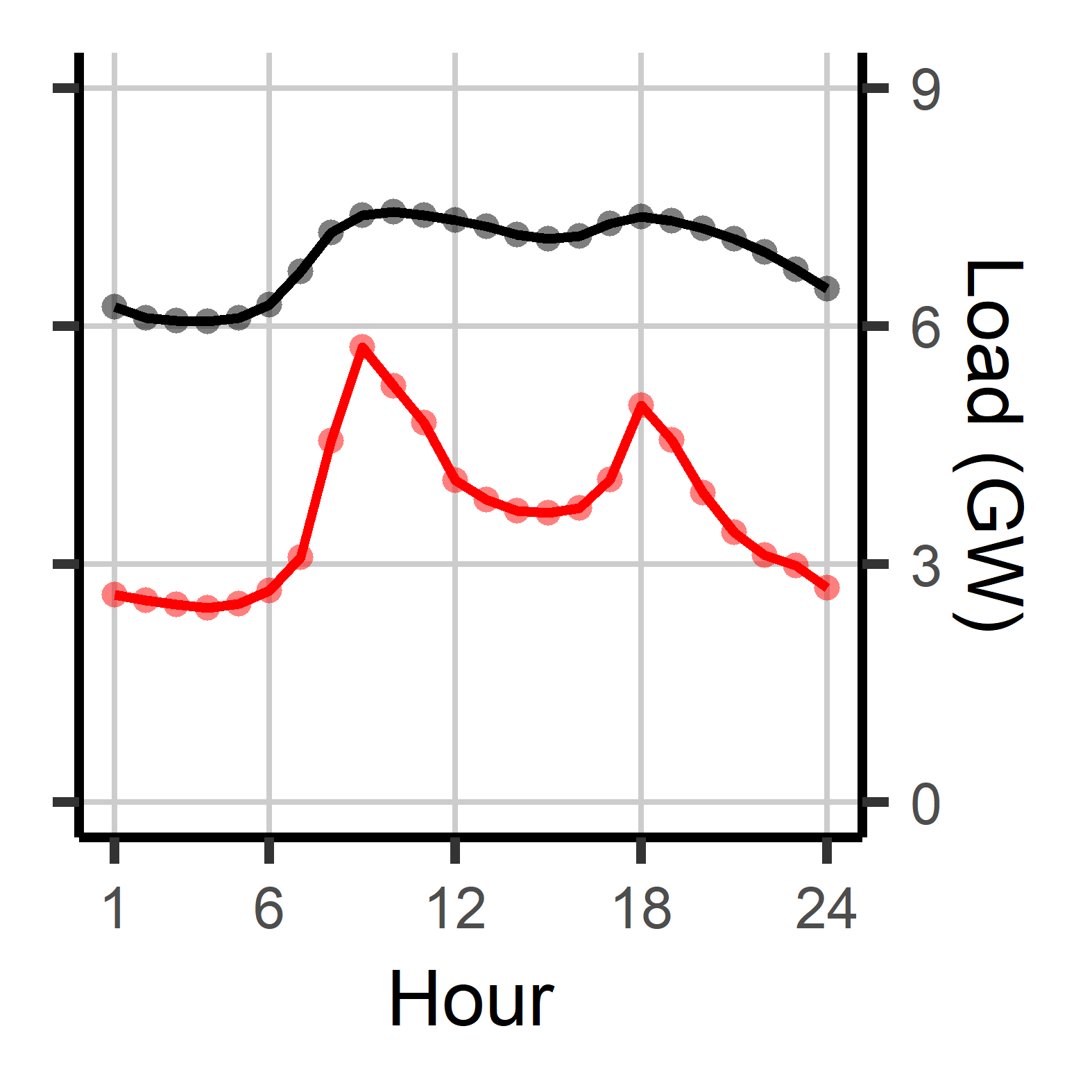}
         \caption{Hourly}
         \label{fig:dailycomparison}
     \end{subfigure}
     \caption{Average daily and hourly load and spot prices in NO1 for a year and all days  with demand over 90\% of peak demand.}
     \label{fig:spot_demand_comparison}
\end{figure}

On peak load days, as shown in \Cref{fig:dailycomparison}, the demand and spot price have the same profile, meaning that the spot price gives the price signal in the right hours. However, when comparing demand and price on the yearly scale, days with high demand does not necessarily imply a high spot price. \Cref{fig:yearlycomparison} shows that the price sometimes correlates well with high demand, but not always. In other words, the spot price might sometimes "interfere" with the grid tariff price signal.

Since the case study looks into daily flexibility, the price signal from spot prices should strengthen the incentive to avoid electricity consumption in peak demand hours on a daily basis. However, the peak demand hours have an almost flat profile, whereas the spot price has larger variations, leading to a strong incentive to reduce demand in a few peak demand hours and not over all peak demand hours.

\subsection{Strength of price signals}
The incentive given through the various grid tariffs or the spot price differs. Based on the specification of the grid tariffs and the historic data, the short-term price signals can theoretically reach values up to 1.2 NOK/kWh for static ToU, 4.5 NOK/kWh for dynamic ToU, 75.25 NOK/kWh/h for demand charges, and 1.65 NOK/kWh for subscribed capacity, compared to the maximum spot price of 2.57 NOK/kWh. On an aggregated level, the price signal is strongest from the dynamic ToU-tariff since the cost are distributed over a few days and are in place for all consumers. In theory, the demand charges gives an even stronger price signal, but only on individual level. Since the monthly peak of all customers is not in the same hour, the aggregated price signal is far lower. 

\begin{figure}[h]
    \centering
    \includegraphics[width=0.8\columnwidth]{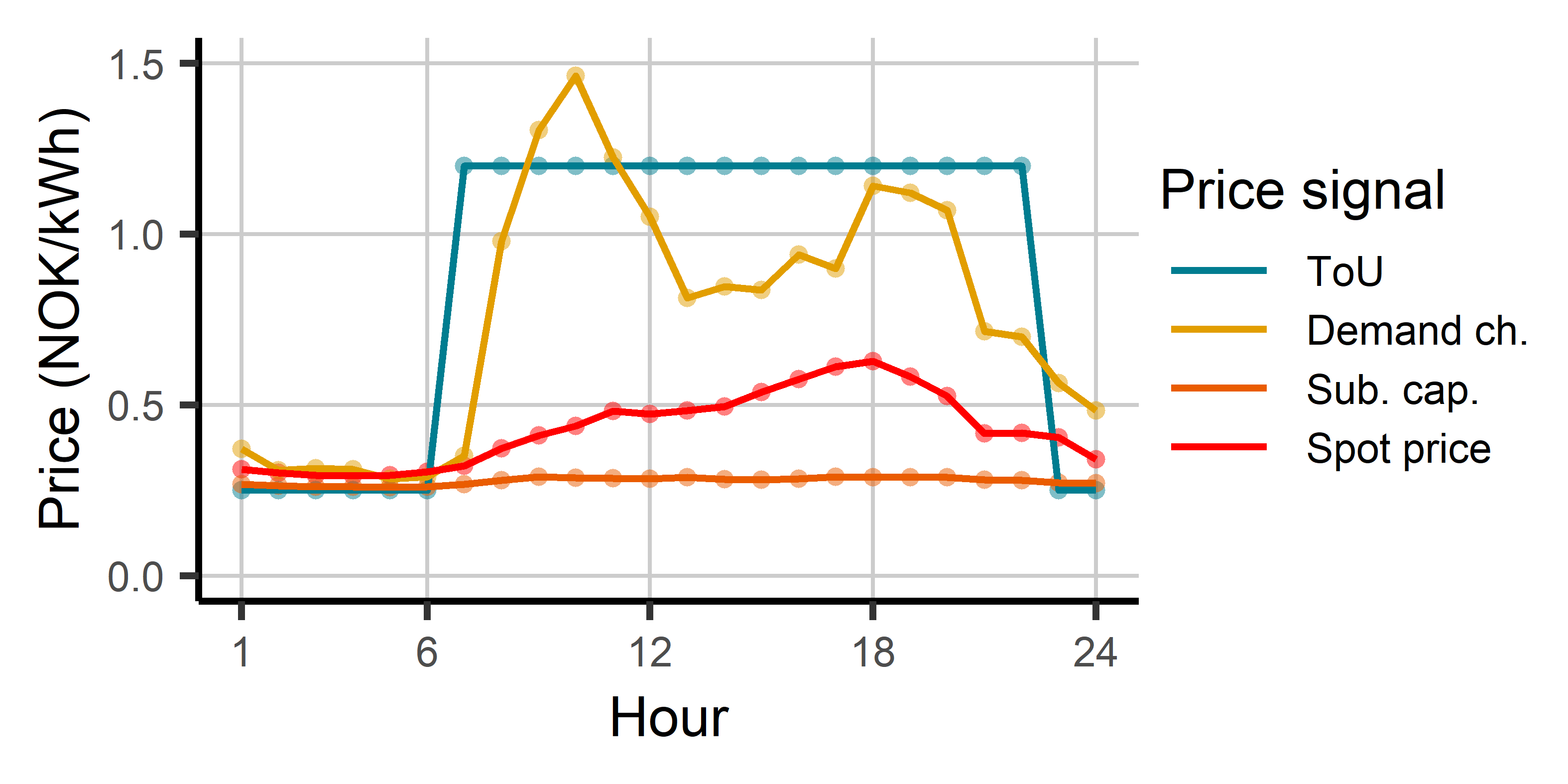}
    \caption{Comparison of price signals from grid tariffs and spot price on the day with maximum demand.}
    \label{fig:prices}
\end{figure}

As an example, Figure \ref{fig:prices} compares the price signals on the aggregated grid level for the day with the highest peak demand for ToU, demand charges, subscribed capacity and spot price. The results show clearly that subscribed capacity gives the weakest price signal on that day, whereas demand charges and the ToU-tariffs give a strong price signal. The spot price is lower than these, but adds an hourly price differentiation to the ToU-tariffs. 

\subsection{Peak change with energy reducing flexibility}
\begin{figure}[h]
    \centering
    \includegraphics[width=0.8\columnwidth]{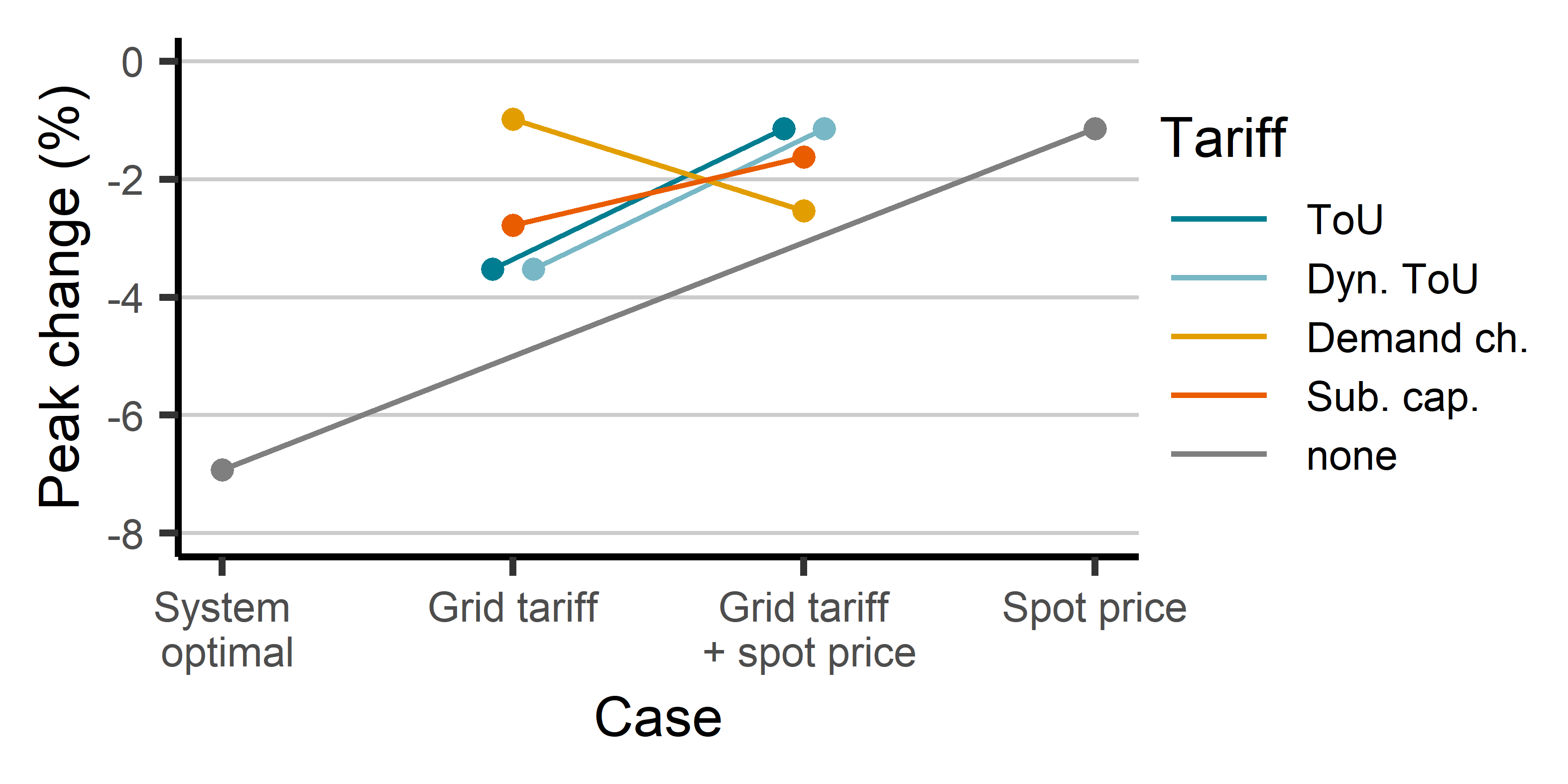}
    \caption{Peak demand reduction results for all cases and grid tariffs.}
    \label{fig:results}
\end{figure}

The results are presented for all cases, i.e. maximum achievable peak reduction with optimal response, price signal from the various grid tariffs, both grid tariffs and spot price, and the spot price alone. The theoretically maximum peak reduction with system optimal response is 6.9 \%, whereas the spot price alone leads to a significant lower reduction of only 1.1 \%. Grid tariffs achieve a peak reduction between 1 to 3.5 \%. However, as Figure \ref{fig:results} shows, all grid tariffs, besides demand charges, achieve a even lower peak reduction when the additional price signal from the spot price is present. The negative effect is largest for the time-of-use tariffs and the achieved peak reduction is then equal to the case with spot price as single price signal. The demand charges tariff leads to the largest peak reduction in combination with the spot price.

The explanation for these results is that time-of-use tariffs have a fixed cost per kWh in the peak hours and therewith, the spot price becomes the predominant price signal in these hours since it varies from hour to hour. In the optimization, all flexibility is therefore used in the hours with the highest spot prices, whereas it otherwise is distributed evenly over all hours with equal peak prices in the ToU-tariffs. Since a reduction of the peak demand needs a load reduction of all hours between 8-21, the achieved peak reduction with spot prices is lower.

\begin{figure}[h]
    \centering
    \includegraphics[]{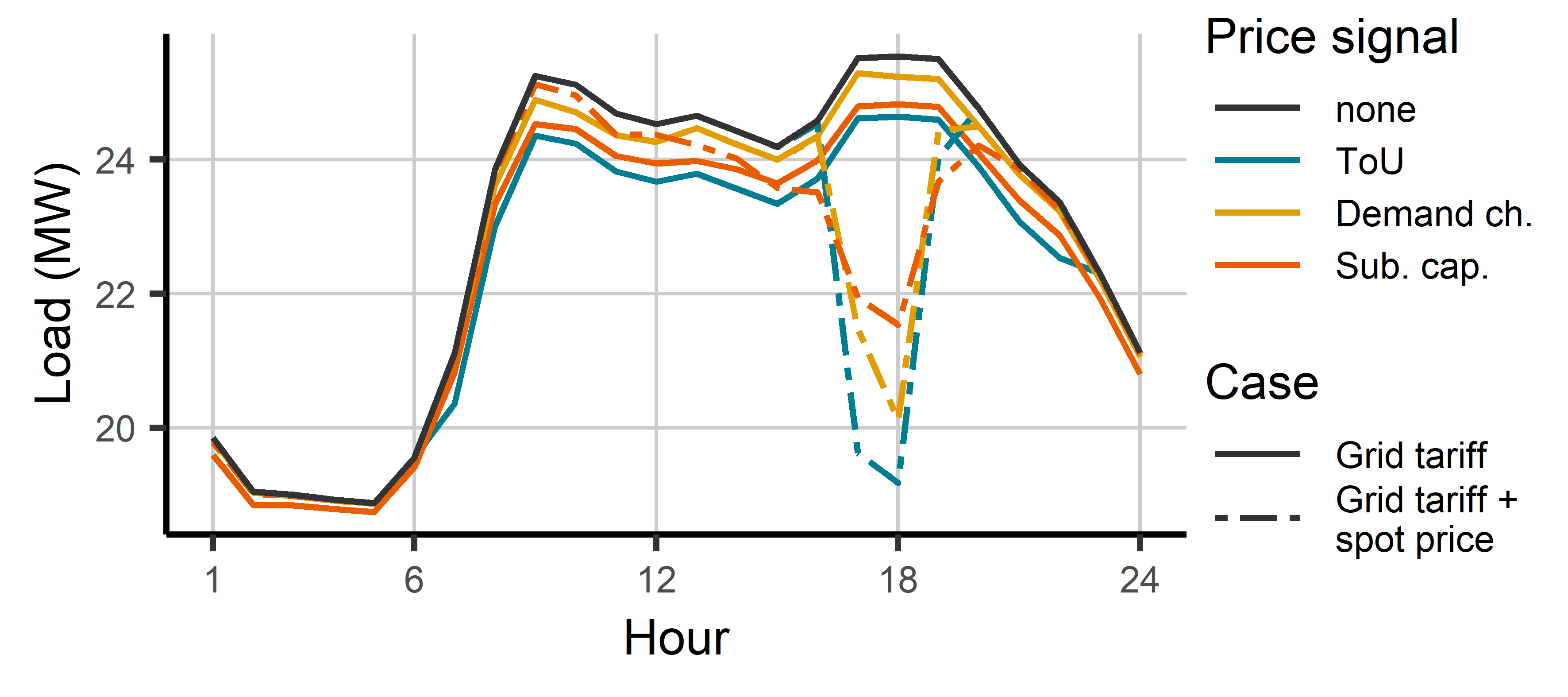}
    \caption{Load change results on the day with maximum demand for various grid tariffs with and without spot price.}
    \label{fig:peakday}
\end{figure}

Figure \ref{fig:peakday} exemplifies these results by showing the load changes on the day with the highest peak demand for the grid tariffs ToU, demand charges and subscribed capacity with and without spot price. The dynamic ToU is not presented since the results are equal to the ToU-tariff. On this day the spot price was highest in the afternoon peak demand hours, but not in the morning peak demand hours. Therefore, in the time-of-use tariffs together with the spot price, the load is only reduced in the afternoon hours. The same effect is also present in the subscribed capacity tariff. However, in the demand charge tariff, only a minor share of the customers have their monthly peak on the grid peak day. Therefore, the additional price signal from the spot price uses mainly unused flexibility to reduce the load in the hours with highest spot price, leading to an increase in peak reduction for this tariff.

\section{Conclusions and further work} \label{sec:conclusion}

This paper demonstrates different grid tariffs designs ability to reduce peak loads, with and without an additional price signal from electricity spot prices. The data shows that in this study, there is no correlation between peak load and peak spot prices over the year, but that the correlation between spot prices and load on peak load days is strong. When subject to electricity spot prices, the consumer demand response leads to smaller reductions in peak load, except for demand charges which performed better together with spot prices. Another conclusion is that even small spot price fluctuations in combination with automatic demand response will lead to use of all flexibility in a few hours. Since the load is high over many hours on peak load days, this leads to inefficient demand response and a low reduction in peak demand. Further work should investigate the impact of different flexibility characteristics from different consumer groups, as this impacts the ability to reduce peak loads.

\section*{Acknowledgements}
This work was funded by the ``DigEco - Digital Economy'' project from the NTNU Digital Transformation Initiative (project number 2495996), Statnett SF, and the Norwegian research council (project number 286513). In addition, it has been supported by CINELDI – Centre for intelligent electricity distribution, an 8-year Research Centre under the FME scheme (Centre for Environment-friendly Energy Research, 257626/E20) of the Norwegian Research Council. The authors thank Kjersti Vøllestad from Elvia for providing the electricity consumption data. 

\bibliographystyle{IEEEtran}
\bibliography{IEEEabrv,mylib}

\appendix

\begin{appendices} \label{Appendix}

\printnomenclature

\end{appendices}

\end{document}